\documentclass{amsart}
\usepackage{graphicx}
\usepackage{amscd}
\usepackage{amsmath}
\usepackage{amsfonts}
\usepackage{amssymb}
\newtheorem{theorem}{Theorem}
\theoremstyle{plain}

\newtheorem{corollary}{Corollary}

\newtheorem{definition}{Definition}
\newtheorem{example}{Example}

\newtheorem{lemma}{Lemma}

\newtheorem{proposition}{Proposition}
\newtheorem{remark}{Remark}

\numberwithin{equation}{section}

\newcommand{\goe}{\mathfrak g}

\newcommand{\ie}{i.e.\ }

\newcommand{\cf}{cf.\ }

\newcommand{\ad}{\operatorname{ad}}

\newcommand{\supp}{\operatorname{supp}}

\begin{document}

\title[Smooth perfectness of diffeomorphisms]{Smooth perfectness through decomposition of diffeomorphisms into fiber
       preserving ones.}
\author{Stefan Haller, Josef Teichmann}

\address{Stefan Haller,
         University of Vienna,
         Strudlhofgasse 4, A-1090 Vienna, Austria.}
\email{stefan@mat.univie.ac.at}
\thanks{The first author is supported by the ``Fonds zur F\"orderung der
        wissenschaftlichen Forschung'' (Austrian Science Fund), project
        number {\tt P14195-MAT}. The second author acknowledges the former supports from the research project {\tt P14195-MAT}
        and from the research grant ``Wittgenstein prize for Prof. Walter Schachermayer'' {\tt Z36-MAT}}

\address{Josef Teichmann,
         Institute of financial and actuarial mathematics,
         Technical University of Vienna,
         Wiedner Hauptstrasse 8--10, A-1040 Vienna, Austria.}
\email{josef.teichmann@fam.tuwien.ac.at}

\begin{abstract}
We show that on a closed smooth manifold $M$ equipped with $k$ fiber
bundle structures whose vertical distributions span the tangent bundle,
every smooth diffeomorphism $f$ of $M$ sufficiently close to the identity can
be written as a product $f=f_1\cdots f_k$, where $f_i$ preserves the
$i^{\text{th}}$-fiber. The factors $f_i$ can be chosen smoothly in $f$.
We apply this result to show that on a certain class of closed smooth
manifolds every diffeomorphism sufficiently close to the identity can be
written as product of commutators and the factors can be chosen smoothly.
Furthermore we get concrete estimates on how many commutators are necessary.
\end{abstract}

\subjclass{58D05}

\maketitle

\section{Introduction}

We are concerned with the question of perfectness of diffeomorphism groups on
compact manifolds. It is well known that the $e$-components of diffeomorphism
groups on compact smooth manifolds are perfect by results of Herman
\cite{Her:73}, Thurston \cite{Thu:74}, Mather \cite{Mat:74,Mat:75}
and Epstein \cite{Eps:84}. However, the questions, how many commutators are
necessary to represent a given smooth diffeomorphism $f$ via
\[
f=[h_{1},g_{1}]\cdots\lbrack h_{n},g_{n}]
\]
and if these commutators can be chosen smoothly in $f$, remains open. Only in
the case of the torus $T^{n}$ the result of Herman provides the concrete,
positive answer by a beautiful small denominator argument.

We shall provide concrete, positive answers for both questions in a
subclass of all compact smooth manifolds by a decomposition theorem
(section~2) and applications of some canonical exponential laws
(``parameterization of diffeomorphisms''). In particular the odd dimensional
spheres $S^{2n+1}$ can be treated, which led us - in view of the applied
methods - to the title of our article.

\section{Inverse Function Theorems}

We shall apply Nash--Moser inverse function theorems in the spirit of Richard
Hamilton, see \cite{Ham:82} for all necessary details, where the theory of
tame spaces and the tame inverse function theorems are presented. Given graded
Fr\'{e}chet spaces $E$ and $F$, \ie we are additionally given an increasing
sequence of seminorms $\{p_{n}\}_{n\geq0}$ on $E$ and $\{q_{n}\}_{n\geq0}$ on
$F$, then a tame estimate for a linear map $L:E\rightarrow F$ is given by%
\[
q_{n}(Le)\leq C_{n}p_{n+r}(e)
\]
for a given basis $b$ and the (tame) degree $r$ for all $n\geq b$ and $e\in
E$. In particular such a map is continuous, we shall call them tame linear
maps. A linear isomorphism $L:E\rightarrow F$ is called a tame isomorphism if
$L$ and $L^{-1}$ satisfy tame estimates. To prove the inverse function theorem
we need to work on tame Fr\'{e}chet spaces: a tame Fr\'echet space is a graded
Fr\'echet space, which is a tame direct summand
of a space $\Sigma(B)$, the Fr\'{e}chet space of all very
fast falling sequences in a Banach space $B$. Let $E$ and $F$ be graded
Fr\'{e}chet spaces and $P:E\supseteq U\rightarrow F$ be a map, then $P$
satisfies a tame estimate if%
\[
q_{n}(P(e))\leq C_{n}(1+p_{n+r}(e))
\]
for a given basis $b$ and the (tame) degree $r$ for all $n\geq b$ and $e\in
U$. We shall call such maps tame maps. Clearly a mapping is tame linear if it
is tame and linear. For mappings on products we can define tame degrees for
any term in the product, which is useful in applications. We shall work in the
category $\mathcal{T}$, whose objects are open subsets of tame spaces and
whose morphisms the smooth tame maps, \ie mappings such that all derivatives
are tame. The Nash--Moser inverse function theorems will be stated in this
category. By means of the tame category
we can define tame manifolds, bundles and geometric or algebraic
structures.

The inverse function theorem finally reads as follows in its general version
in $\mathcal{T}$, see \cite{Ham:82} for both theorems and all necessary
details.

\begin{theorem}
Let $E$ and $F$ be tame spaces and $P:E\supseteq U\rightarrow F$ be a smooth
tame map. Suppose that the equation $DP(e)h=f$ has a unique solution
$h=VP(e)f\in E$ for all $e\in U$ and all $f\in F$ and that the family of inverses
$VP:U\times F\rightarrow E$ is smooth tame, then $P$ is locally invertible and
the inverse is a smooth tame map.
\end{theorem}

We shall apply the version for right inverses.

\begin{theorem}
Let $E$ and $F$ be tame spaces and $P:E\supseteq U\rightarrow F$ be a smooth
tame map. Suppose that the equation $DP(e)h=f$ has a solution $h=VP(e)f\in E$ for
all $e\in U$ and all $f\in F$ and that the family of right inverses
$VP:U\times F\rightarrow E$ is smooth tame, then $P$ is locally surjective and
admits a smooth tame local right inverse.
\end{theorem}

\section{The Decomposition Theorem}

In this section we shall prove the fundamental decomposition theorem, which
allows a ``smooth'' decomposition of a small diffeomorphism on a compact
manifold into more regular parts. This shall be applied in the following
sections to obtain perfectness results on Fr\'{e}chet--Lie groups such as
$\text{\textrm{Diff}}(S^{2n+1})$.

Let $M$ be a closed connected manifold such that there exist $k$ fiber bundle
structures $S_{i}\hookrightarrow M\overset{p_{i}}{\longrightarrow}B_{i}$ for
$1\leq i\leq k$ with connected fibers $S_{i}$ and (involutive) vertical
distribution $\mathcal{D}_{i}$, which is a subbundle of $TM$.
We suppose that the distributions $\mathcal{D}_{i}$ span $TM$, however they
need not be linearly independent.\footnote{Note,
that if $M$ appears as the total space of a fiber bundle
$S\hookrightarrow M\overset{p}\longrightarrow B$ and $\dim S\geq 1$,
then one can always perturb $p$ to obtain finitely many fiber bundles
$S\hookrightarrow M\overset{p_i}\longrightarrow B$, which satisfy this
condition.} We denote the Lie subgroup of $\operatorname{Diff}(M)$
of bundle diffeomorphisms by $\operatorname{Diff}_{i}(M)$ for $1\leq i\leq k$.
The Lie algebra of $\operatorname{Diff}(M)$ is given by smooth sections
$\Gamma(TM)$, the Lie algebra of $\operatorname{Diff}_{i}(M)$
by sections $\Gamma(\mathcal{D}_{i})$. These Lie groups are
tame manifolds, \ie a smooth tame atlas exists and the structure maps are
smooth tame. Remark in particular that the pullback, \ie the adjoint action
$\operatorname{Ad}:\operatorname{Diff}(M)\times\Gamma(TM)\to\Gamma(TM)$,
$(f,X)\mapsto (f^{-1})^*X$
is smooth tame as derivative of the conjugation. Furthermore the
module structure on $\Gamma(TM)$ is tame. See \cite{KriMic:97}
for the general theory of Lie groups and \cite{Ham:82} for tame Lie groups.

We shall apply the Nash--Moser Theorem in the following version: Given
tame manifolds $\mathcal{M}$ and $\mathcal{N}$ and a smooth tame map
$P:\mathcal{M}\rightarrow\mathcal{N}$, then the existence of a local
smooth tame right inverse of $P$ is equivalent to the existence of a
local smooth tame right inverse (vector bundle map) to
$\widetilde{TP}:T\mathcal{M}\rightarrow P^{\ast}T\mathcal{N}$, where
$T\mathcal{N}$ denotes the (tame) tangent bundle and
$P^{\ast}T\mathcal{N}$ denotes the (tame) pullback bundle.

\begin{theorem}\label{decomposition_theorem}
The smooth tame mapping%
\begin{gather*}
P:\operatorname{Diff}_{1}(M)\times\cdots\times\operatorname{Diff}_{k}
(M)\rightarrow\operatorname{Diff}(M)
\\
(f_{1},\dotsc,f_{k})\mapsto f_{1}\circ\cdots\circ f_{k}
\end{gather*}
admits a smooth tame local right inverse at the identity
$e\in\operatorname{Diff}(M)$.
\end{theorem}

\begin{proof}
In the right trivialization of the tangent bundles of the respective Lie
groups we are given a mapping:
\begin{gather*}
\widetilde{TP}:
\prod_{i=1}^{k}\operatorname{Diff}_{i}(M)\times
\prod_{i=1}^{k}\Gamma(\mathcal{D}_{i})
\rightarrow
\prod_{i=1}^{k}\operatorname{Diff}_{i}(M)\times\Gamma(TM)
\\
(f_{1},\dotsc,f_{k};\xi_{1},\dotsc,\xi_{k})
\mapsto
(f_{1},\dotsc,f_{k};f_{2}^{\ast}\cdots f_{k}^{\ast}\xi_{1}+
f_{3}^{\ast}\cdots f_{k}^{\ast}\xi_{2}+\cdots+\xi_{k})
\end{gather*}
In view of the implicit function theorem it suffices to construct a
smooth tame right inverses of $\widetilde{TP}$, linear in the variables
$\xi_i$. We solve the problem locally: First we choose a covering
$\mathfrak{U}$ of open subsets of $M$ such that for $U\in\mathfrak{U}$
the bundles $TM|_{U}$ and $\mathcal{D}_{i}|_{U}$ are trivial and that a
compatible local frame (see below) exists. Second we choose a finite
partition of unity $\{\eta_{U}\}_{U\in\mathfrak{U}}$ subordinated to
this covering. The
associated projection%
\begin{gather*}
p:\Gamma(TM)\rightarrow\bigoplus_{U\in\mathfrak{U}}\Gamma_{V}(TM|_{U})
\\
p_{U}(X)=\eta_{U}X
\end{gather*}
where $\Gamma_{V}$ denotes the sections with support in the closed set $V$ and
$V$ is the support of $\eta_{U}$ in $U$, is a right inverse of the sum.
On $U\in\mathfrak{U}$
we can solve the equation locally. We choose a local frame $X^{1},\dotsc
,X^{n}$, where $n=\dim M$, compatible with the distribution on $U$, \ie
there are integers $0=m_{0}\leq m_{1}\leq\cdots\leq m_{k}=n$, such that
\[
\mathcal{D}_{i}(x)
=\left\langle X^{n_{i}}(x),\dotsc,X^{m_{i}}(x)\right\rangle
\quad\text{for all $x\in U$,}
\]
where we set $n_i:=m_{i-1}+1$. We
assume furthermore that the vector fields $X^{j}$ are globally defined on $M$.
We then choose an open neighborhood $V_{i}$ of $e\in\operatorname{Diff}%
_{i}(M)$ such that for all $f_{i}\in V_{i}$ the condition $f_{1}\circ
\cdots\circ f_{n}(V)\subset U$ holds and such that%
\[
\mathfrak{F}_{(f_{1},\dotsc,f_{k})}
:=(f_{2}^{\ast}\cdots f_{k}^{\ast}X_{n_{1}},\dotsc,
f_{2}^{\ast}\cdots f_{k}^{\ast}X_{m_{1}},
f_{3}^{\ast}\cdots f_{k}^{\ast}X_{n_{2}},\dotsc,X_{n})
\]
is a frame for $TM$ on $V$. Given $Y\in\Gamma_{V}(TM|_{U})$,
we define a section $s_{i,U}$ via the decomposition on the frame $\mathfrak
{F}_{(f_{1},\dotsc,f_{k})}$ by the following formula%
\[
s_{i,U}(Y):=(f_{i+1}^{\ast}\cdots f_{k}^{\ast})^{-1}
\left(
\sum_{j=n_{i}}^{m_{i}}
a_{j}(Y)f_{i+1}^{\ast}\cdots f_{k}^{\ast}X_{j}
\right)
\]
where
\[
Y=\sum_{i=1}^{k}\sum_{j=n_{i}}^{m_{i}}a_{j}(Y)f_{i+1}^{\ast}\cdots
f_{k}^{\ast}X_{j}\text{.}%
\]
The functions $a_{j}$ are smooth with support in $V$, so $s_{i,U}(Y)$
has support in $U$ and consequently the section defines an element of
$\Gamma(\mathcal{D}_{i})$ by construction. The pullbacks depend tame on
the diffeomorphism and the module structure is tame, hence the
functions $a_{j}$ depend tame on the diffeomorphisms
$f_{i+1},\dotsc,f_{k}$ and the vector field $Y$. Consequently the
mapping
\begin{gather*}
\prod_{i=1}^{k}V_i\times\Gamma(TM)
\rightarrow
\prod_{i=1}^{k}\operatorname{Diff}_{i}(M)\times
\prod_{i=1}^{k}\Gamma(\mathcal{D}_{i})
\\
(f_{1},\dotsc,f_{k};X)
\mapsto
(f_{1},\dotsc,f_{k};\sum_{U\in\mathfrak{U}}s_{1,U}(\eta_{U}X),\dotsc,
\sum_{U\in\mathfrak{U}}s_{k,U}(\eta_{U}X))
\end{gather*}
is the desired smooth tame right inverse of $\widetilde{TP}$.
\end{proof}

\section{Perfectness of fiber-preserving Diffeomorphisms}

In this section we introduce notions of perfectness on regular Lie
groups, namely global, semi infinitesimal and infinitesimal smooth
perfectness. The relations between these notions become complicated on the
level of Fr\'{e}chet--Lie groups, even though natural inequalities remain
valid. These notions are finally applied to ``parametrized'' families of
diffeomorphisms, which is an application of cartesian closedness, see for
example \cite{KriMic:97}.

All manifolds, Lie groups and Lie algebras are supposed to be smooth and
modeled on convenient vector spaces. This includes all Fr\'echet manifolds,
Fr\'echet--Lie groups and Fr\'echet--Lie algebras. We moreover assume, that
the Lie groups are regular, in particular they admit a smooth exponential
mapping. This is not too much a restriction, since
``all known convenient Lie groups are regular'', \cf \cite{KriMic:97}.
Conditions for regularity of convenient Lie groups can be found in
\cite{Omo:97}, \cite{Tei:01}.

\begin{definition}
For a Lie group $G$ we define $N_{G}\in\mathbb{N}$ to be the smallest integer
$N$, such that for every open neighborhood $e\in U\subseteq G$ their exist
$h_{i}=\exp(Y_i)\in U$, an open neighborhood $e\in V\subseteq G$ and smooth
mappings
$S_{i}:V\rightarrow G$ with $S_{i}(e)=e$ and
\[
\lbrack S_{1}(g),h_{1}]\cdots\lbrack S_{N}(g),h_{N}]=g,\quad\text{for all
$g\in V$.}%
\]
Equivalently, $N_{G}$ is the smallest integer $N$, such that for every open
neighborhood $e\in U\subseteq G$ there exist $h=\exp(Y)\in U^{N}$ with, such
that the map
\[
\kappa_{h}:G^{N}\rightarrow G,\quad(g_{1},\dotsc,g_{N})\mapsto\lbrack
g_{1},h_{1}]\cdots\lbrack g_{N},h_{N}]
\]
has a smooth local right inverse $S$ with $S(e)=(e,\dotsc,e)$. If such an
integer does not exist we set $N_{G}:=\infty$. We call the Lie group $G$
smoothly perfect if $N_{G}<\infty$.
\end{definition}

\begin{definition}
For a Lie group $G$ with Lie algebra $\mathfrak g$ we define $N_{G}
^{\operatorname{Ad}}\in\mathbb{N}$ to be the smallest integer $N$, such that
for every open neighborhood $e\in U\subseteq G$ their exist $h_{i}=\exp
(Y_{i})\in U$ and bounded linear maps $s_{i}:\mathfrak g\rightarrow\mathfrak
g$ with
\[
(\operatorname{id}-\operatorname{Ad}_{h_{1}})s_{1}(X)+\cdots
+(\operatorname{id}-\operatorname{Ad}_{h_{N}})s_{N}(X)=X\quad\text{for all
$X\in\mathfrak g$.}%
\]
Equivalently, $N_{G}^{\operatorname{Ad}}$ is the smallest integer $N$, such
that for every open neighborhood $e\in U\subseteq G$ there exist
$h=\exp(Y)\in U^{N}$ and a bounded linear right inverse $s:\mathfrak
g\rightarrow\mathfrak g^{N}$ of the map $T_{(e,\dotsc,e)}\kappa_{h}:\mathfrak
g^{N}\rightarrow\mathfrak g$. If such an integer does not exist we set
$N_{G}^{\operatorname{Ad}}:=\infty$.
\end{definition}

\begin{remark}
It would be more natural to claim existence of arbitrary, small $h_{i}$, not
only those which are exponentials, however, in this general case we do not get
the desired ``natural'' inequalities, \cf Lemma~\ref{NgoevsNG} below.
\end{remark}

\begin{definition}
For a Lie algebra $\mathfrak g$ we define $N_{\mathfrak g}\in{\mathbb N}$
to be the
smallest integer $N$, such that there exist $Y_{i}\in\mathfrak g$ and bounded
linear maps $\sigma_{i}:\mathfrak g\rightarrow\mathfrak g$ with
\[
\lbrack\sigma_{1}(X),Y_{1}]+\cdots+[\sigma_{N}(X),Y_{N}]=X,\quad\text{for all
$X\in\mathfrak g$.}%
\]
Equivalently $N_{\mathfrak g}$ is the smallest integer $N$, such that there
exist $Y\in\mathfrak g^{N}$ and a bounded linear right inverse $\sigma
:\mathfrak
g\rightarrow\mathfrak g^{N}$ of the mapping
\[
K_{Y}:\mathfrak g^{N}\rightarrow\mathfrak g,\quad K_{Y}(X_{1},\dotsc
,X_{N}):=[X_{1},Y_{1}]+\cdots+[X_{N},Y_{N}].
\]
If such an integer does not exist we set $N_{\mathfrak g}:=\infty$.
\end{definition}

\begin{lemma}\label{NgoevsNG}
For any Lie group $G$ with Lie algebra $\mathfrak g$ one has $N_{\mathfrak
g}\leq N_{G}^{\operatorname{Ad}}\leq N_{G}$. If $G$ is a Banach--Lie group one
even has $N_{\mathfrak g}=N_{G}^{\operatorname{Ad}}=N_{G}$.
\end{lemma}

\begin{proof}
$N_{G}^{\operatorname{Ad}}\leq N_{G}$ follows immediately from differentiating
$\kappa_{h}\circ S=\text{\textrm{id}}$ at $e\in G$, i.e.\ one can take
$s=T_{e}S$. If $G$ is a Banach--Lie group then the implicit function theorem shows
$N_{G}^{\operatorname{Ad}}\geq N_{G}$. Notice, that every $h$ sufficiently
close to the identity is in the image of the exponential map, for the latter
is a local diffeomorphism on Banach--Lie groups.

Next we show $N_{\mathfrak g}\leq N_{G}^{\operatorname{Ad}}$. For
$Y\in{\mathfrak g}$ we have
$\frac{\partial}{\partial t}\operatorname{Ad}_{\exp(tY)}
=\operatorname*{ad}_{Y}\circ\operatorname{Ad}_{\exp(tY)}$.
Integration immediately yields
\[
\operatorname{id}-\operatorname{Ad}_{h}
= - {\operatorname*{ad}}_{Y} \circ \int_{0}^{1}%
\operatorname{Ad}_{\exp(tY)}dt,
\]
where $h=\exp(Y)$. Inserting the bounded linear right inverse $s$ for
$T_{(e,\dotsc,e)}\kappa_h$,
we obtain a bounded linear right inverse
$\sigma=(\int_{0}^{1}\operatorname{Ad}_{\exp(tY)}dt\circ s_{i})_{i=1,\dotsc,N_{\mathfrak g}}$ for
$K_{Y}$.

Suppose $G$ is a Banach--Lie group. We want to show $N_{\mathfrak g}\geq
N_{G}^{\operatorname{Ad}}$. Choose $Y\in\mathfrak g^{N_{\mathfrak g}}$, such
that $K_{Y}$ has a bounded linear right inverse and choose a smooth curve
$h_{t}\in U^{N_{\mathfrak g}}$ with $h_{0}=(e,\dotsc,e)$ and $\dot{h}_{0}=Y$.
For $t>0$ consider the maps $K_{t}:=\frac{1}{t}T_{(e,\dotsc,e)}\kappa_{h_{t}}$
and note, that $\lim_{t\rightarrow0}K_{t}=K_{Y}$. Since $\mathfrak g$ is a
Banach space, the space of bounded linear mappings $\mathfrak g^{N_{\mathfrak
g}}\rightarrow\mathfrak g$ which admit a bounded linear right inverse is open,
hence for $t$ sufficiently small $K_{t}$ has a bounded linear right inverse
and thus $T_{(e,\dotsc,e)}\kappa_{h_{t}}=tK_{t}$ as well, i.e.\ $N_{\mathfrak
g}\geq N_{G}^{\operatorname{Ad}}$.
\end{proof}

\begin{example}
For a finite dimensional perfect Lie group $G$ one has $1<N_G\leq\dim G$.
For any finite dimensional Lie algebra one has $N_{\mathfrak g}>1$, since
$\ad_Y:{\mathfrak g}\to{\mathfrak g}$ can't be
surjective, for it has a non-trivial kernel.
Moreover obviously $N_{\mathfrak{g}}\leq\dim\mathfrak{g}$.
\end{example}

\begin{example}\label{semisimple}
If $G$ is complex semisimple or real semisimple and split or real
semisimple and compact then $N_{G}=2$. Indeed, if $G$ is complex for
example, choose $H$ to be a regular element in the Cartan subalgebra
${\mathfrak h}$ and $\rho:=\sum_{\alpha}E_\alpha$, where the sum is over all
simple roots and $E_\alpha$ denotes a non-zero element of the root space of
$\alpha$. Then $\operatorname{ad}_H({\mathfrak g})={\mathfrak h}^\perp$
and $\operatorname{ad}_\rho(\goe)\supseteq{\mathfrak h}$, hence
$\operatorname{ad}_H+\operatorname{ad}_\rho:
{\mathfrak g}\times{\mathfrak g}\to{\mathfrak g}$ is onto.
In the two real cases one can argue similarly.
\end{example}

The first non-trivial example is an immediate consequence of a theorem
due to Herman, \cf \cite{Her:73}.

\begin{example}\label{Herman}
For the torus $T^n$ one has $N_{\operatorname{Diff}(T^n)}\leq 3$.
Indeed, Herman proves the statement with one commutator up to multiplication
by an element in $T^n$. Since $T^n\subseteq\operatorname{PSL}(2,{\mathbb R})^n$
and $\operatorname{PSL}(2,{\mathbb R})^n$ is real semisimple and split
example~\ref{semisimple} shows that $N_{\operatorname{Diff}(T^n)}\leq 3$.
\end{example}

The base of all bundles we shall consider below is understood to be a compact,
smooth and finite dimensional manifold, but the fiber might be infinite dimensional.

\begin{definition}
Let $\pi:E\rightarrow B$ be a fiber bundle with typical fiber $F$, whose
structure group is reduced to $K\subseteq\operatorname{Diff} (F)$, i.e.\ we
have given a fiber bundle atlas whose transition functions take values in
$K\subseteq\operatorname{Diff} (F)$, where $K$ is any subgroup of
$\operatorname{Diff} (F)$. We define $C_{E}=C_{E}^{K}$ to be the smallest
integer, such that there exists an open covering $\{U_{1},\dotsc,U_{C_{E}}\}$
of $B$ and a fiber bundle atlas $\varphi_{i}:E|_{U_{i}}\rightarrow
U_{i}\times F$ whose transition functions take values in $K\subseteq
\operatorname{Diff}(F)$.
\end{definition}

\begin{remark}
Note, that we do not assume the $U_{i}$ to be connected. Since every
manifold $B$ can be covered by $\dim B+1$ open sets each of which is a
disjoint union of disks one gets $C_E\leq\dim B+1$ for any bundle
$E\to B$.
\end{remark}

Suppose we have a bundle of Lie groups $E\to B$ with typical fiber $G$,
i.e.\ the structure group is reduced to $\operatorname{Aut}(G)\subseteq
\text{\textrm{Diff}} (G)$. Then the space of smooth sections $\Gamma(E)$ is
again a manifold, which becomes a Lie group under point wise multiplication.

\begin{proposition}\label{prop}
Suppose $\pi:E\rightarrow B$ is a bundle of Lie groups with typical fiber $G$.
Then $N_{\Gamma(E)}\leq C_{E}^{\operatorname{Aut}(G)}N_{G}$.
\end{proposition}

The proposition will follow immediately from the following two lemmas.

\begin{lemma}\label{first_lemma}
Let $E\to B$ be a bundle of Lie groups with typical fiber $G$ and suppose
$\{V_{1},\dotsc,V_{N}\}$ is an open covering of $B$, such that $E|_{V_i}$ is
trivial. Then there exist an open
neighborhood $e\in\mathcal{V}\subseteq\Gamma(E)$ and smooth mappings
$F_{i}:\mathcal{V}\to\Gamma_{\bar V_{i}}(E)$ with $F_{i}(e)=e$ and
$F_{1}(s)\cdots F_{N}(s)=s$, for all $s\in\mathcal{V}$.
\end{lemma}

\begin{lemma}\label{second_lemma}
Suppose $W$ is a finite dimensional manifold which need not be compact, $G$ a
Lie group, $N_{G}<\infty$, $V\subseteq U\subseteq W$ open, such that $\bar
V\subseteq U$ and such that $\bar U$ is compact. Then for every open
neighborhood of $e\in\mathcal{U}\subseteq C^{\infty}_{\bar U}(W,G)$ there
exist $h^{j}=\exp(Y^j)\in\mathcal{U}$, an open neighborhood $e\in\mathcal{V}\subseteq
C^{\infty}_{\bar V}(W,G)$ and smooth mappings $S^{j}:\mathcal{V}\to C^{\infty
}_{\bar V}(W,G)$ with $S^{j}(e)=e$ and $[S^{1}(f),h^{1}]\cdots[S^{N_{G}%
}(f),h^{N_{G}}]=f$, for all $f\in\mathcal{V}$.
\end{lemma}

\begin{proof}[Proof of Lemma~\ref{first_lemma}]
Choose bump functions $\lambda_i:B\to[0,1]$ with
$\supp\lambda_i\subseteq V_i$ and such that $U_i:=\{x\in
B:\lambda_i(x)=1\}$ still cover $B$. Using trivializations of
$E|_{V_i}$ and a chart of $G$ centered at $e$ which has a convex image
in $\mathfrak g$ one defines a smooth map given by ``multiplication
with $\lambda_i$'' $$ \phi_i:\Gamma(E)\supseteq{\mathcal
V}\to\Gamma_{\bar V_i}(E),\quad 1\leq i\leq N, $$ where $\mathcal V$ is
an open neighborhood of the identical section. Obviously that map has
the property, that $\phi_i(s)=s$ on $U_i$ and
$\supp(\phi_i(s))\subseteq\supp(s)$. Now set $F_1(s):=\phi_1(s)$ and $$
F_i(s):=\phi_i(F_{i-1}(s)^{-1}\cdots F_1(s)^{-1}s),\quad 1\leq i\leq N.
$$ Shrinking $\mathcal V$ we may assume that everything is well
defined. An easy inductive argument shows $F_1(s)\cdots F_i(s)=s$ on
$U_1\cup\cdots\cup U_i$, for all $s\in\mathcal V$ and all $1\leq i\leq
N$.
\end{proof}

\begin{proof}[Proof of Lemma~\ref{second_lemma}]
Choose a bump function $\mu:W\to$ with
$\operatorname{supp}(\mu)\subseteq U$ and $\mu=1$ on $V$. Let $\tilde
V$, $\tilde h_{i}=\exp(\tilde Y_i)$ and $\tilde S_{i}:\tilde V\to\tilde U$, $1\leq i\leq N_{G}%
$, be the data we get from $N_{G}<\infty$ and a sufficiently small
neighborhood of $e\in G$. Set $\mathcal{V}%
:=\{f\in C^{\infty}_{\bar V}(W,G): f(\bar V)\subseteq\tilde V\}$,
$h^{i}(x):=\exp(\mu(x)\tilde Y_i)$ and $S^{i}:=(\tilde S_{i})_{*}$.
\end{proof}

\begin{proof}[Proof of Proposition~\ref{prop}]
Choose open sets $V_{i}\subseteq\bar{V}_{i}\subseteq U_{i}\subseteq\bar{U}%
_{i}\subseteq W_{i}$, $1\leq i\leq N_{E}$, such that $E|_{W_{i}}$ are trivial
and such that $\{V_{i}\}$ is an open covering of $B$. Suppose we have given
any open neighborhood $e\in\mathcal{U}\subseteq\Gamma(E)$. For every $1\leq
i\leq N_{E}$ the second lemma provides, via trivializations of $E|_{W_{i}}$,
$h_{i}^{j}\in\mathcal{U}$, an open neighborhood $e\in\mathcal{V}_{i}%
\subseteq\Gamma_{\bar{V}_{i}}(E)$ and smooth mappings $S_{i}^{j}%
:\mathcal{V}_{i}\rightarrow\Gamma(E)$, $1\leq j\leq N_{G}$, with $S_{i}%
^{j}(e)=e$ and $[S_{i}^{1}(s),h_{i}^{1}]\cdots\lbrack S_{i}^{N_{G}}%
(s),h_{i}^{N_{G}}]=s$, for all $s\in\mathcal{V}_{i}$. Let $e\in\mathcal{V}%
\subseteq\Gamma(E)$ be the open neighborhood from the first lemma and assume,
that $F_{i}:\mathcal{V}\rightarrow\mathcal{V}_{i}$. Then $S_{i}^{j}\circ
F_{i}:\mathcal{V}\rightarrow\Gamma(E)$, $(S_{i}^{j}\circ F_{i})(e)=e$ and
\[
\prod_{1\leq i\leq N_{E}}\,\prod_{1\leq j\leq N_{G}}[(S_{i}^{j}\circ
F_{i})(s),h_{i}^{j}]=s,
\]
for all $s\in\mathcal{V}$.
\end{proof}

\begin{example}
Suppose $M\to B$ is a finite dimensional principal bundle with perfect
structure group $G$. Let $\operatorname{Gau}(M)$ denote the group of
gauge transformations. It is well known the
$\operatorname{Gau}(M)=\Gamma(E)$, where $E$ is the associated bundle
of groups with typical fiber $G$. So Proposition~\ref{prop} implies,
that $\operatorname{Gau}(M)$ is a smoothly perfect group, even gives
concrete estimates, e.g.\ at least
$N_{\operatorname{Gau}(M)}\leq(\dim B+1)(\dim G)$.
\end{example}

Similarly one can treat other examples, such as the group of
automorphisms on a finite dimensional vector bundle preserving a fiber metric
(fiber volume, fiber symplectic form), or the group of automorphisms on a
finite dimensional bundle of groups (those which are fiber wise group
isomorphisms), or the group of automorphisms on a finite dimensional bundle
of Lie algebras.

\section{Applications}

Applications are given by a combination of the decomposition theorem and the
fact that%
\[
\operatorname{Diff}(M,S)=\Gamma(E),
\]
where $M$ is a compact fiber bundle $S\hookrightarrow M\rightarrow B$,
$\operatorname{Diff}(M,S)$ denotes the fiber preserving diffeomorphisms of $M$
and $E$ denotes the associated bundle of Lie groups over $B$ with typical
fiber $\operatorname{Diff}(S)$. So Proposition~\ref{prop} shows, that
$\operatorname{Diff}(M,S)$ is smoothly perfect as soon as
$\operatorname{Diff}(S)$ is smoothly perfect. See also \cite{Ryb:95}, where
it is shown, that the group of leaf preserving diffeomorphisms on any foliated
manifold is perfect. From Theorem~\ref{decomposition_theorem} and
Proposition~\ref{prop} we immediately get the following

\begin{corollary}
Suppose $M$ is a closed manifold which admits $k$ fiber bundles
$S_i\hookrightarrow M\overset{p_i}\longrightarrow B_i$ such that the corresponding vertical distributions span
$TM$. Then
$$
N_{\operatorname{Diff}(M)}
\leq
\sum_{i=1}^kC_{p_i}N_{\operatorname{Diff}(S_i)}.
$$
Particularly $\operatorname{Diff}(M)$ is smoothly perfect if all
$\operatorname{Diff}(S_i)$ are smoothly perfect.
\end{corollary}

Note, that if $M$ appears as a total space of one fiber bundle
$S\hookrightarrow M\overset{p}\longrightarrow B$ with $\dim(S)\geq 1$,
one can always perturb
$p$ and find finitely many fiber bundles $S\hookrightarrow
M\overset{p_i}\longrightarrow B$, whose vertical distributions span $TM$.
Moreover we always have the estimate $C_p\leq\dim B+1\leq\dim M$.

\begin{example}
Since every odd dimensional sphere appears as the total space of an
$S^1$-bundle via the Hopf fibration Example~\ref{Herman} implies, that
$\operatorname{Diff}(S^{2n+1})$ is smoothly perfect. For example one
easily derives $N_{\operatorname{Diff}(S^3)}\leq 18$.
\end{example}

\begin{example}
Since every compact Lie group $G$ has a torus as subgroup it appears as total
space in an $S^1$-bundle and hence $\operatorname{Diff}(G)$ is smoothly
perfect. We even get the estimate
$N_{\operatorname{Diff}(G)}\leq 3(\dim G)^2$, for there are always $\dim G$
many $S^1$-bundle structures on $G$ which span $TG$.
\end{example}

\end{document}